\documentclass{amsart}
\usepackage{amssymb}

\makeatletter
 
\def\diagram{\m@th\leftwidth=\z@ \rightwidth=\z@ \topheight=\z@
\botheight=\z@ \setbox\@picbox\hbox\bgroup}
 
\def\enddiagram{\egroup\wd\@picbox\rightwidth\unitlength
\ht\@picbox\topheight\unitlength \dp\@picbox\botheight\unitlength
\hskip\leftwidth\unitlength\box\@picbox}
 
\def\bfig{\begin{diagram}}
\def\efig{\end{diagram}}
\newcount\wideness \newcount\leftwidth \newcount\rightwidth
\newcount\highness \newcount\topheight \newcount\botheight
 
\def\ratchet#1#2{\ifnum#1<#2 \global #1=#2 \fi}
 
\def\putbox(#1,#2)#3{%
\horsize{\wideness}{#3} \divide\wideness by 2
{\advance\wideness by #1 \ratchet{\rightwidth}{\wideness}}
{\advance\wideness by -#1 \ratchet{\leftwidth}{\wideness}}
\vertsize{\highness}{#3} \divide\highness by 2
{\advance\highness by #2 \ratchet{\topheight}{\highness}}
{\advance\highness by -#2 \ratchet{\botheight}{\highness}}
\put(#1,#2){\makebox(0,0){$#3$}}}
 
\def\putlbox(#1,#2)#3{%
\horsize{\wideness}{#3}
{\advance\wideness by #1 \ratchet{\rightwidth}{\wideness}}
{\ratchet{\leftwidth}{-#1}}
\vertsize{\highness}{#3} \divide\highness by 2
{\advance\highness by #2 \ratchet{\topheight}{\highness}}
{\advance\highness by -#2 \ratchet{\botheight}{\highness}}
\put(#1,#2){\makebox(0,0)[l]{$#3$}}}
 
\def\putrbox(#1,#2)#3{%
\horsize{\wideness}{#3}
{\ratchet{\rightwidth}{#1}}
{\advance\wideness by -#1 \ratchet{\leftwidth}{\wideness}}
\vertsize{\highness}{#3} \divide\highness by 2
{\advance\highness by #2 \ratchet{\topheight}{\highness}}
{\advance\highness by -#2 \ratchet{\botheight}{\highness}}
\put(#1,#2){\makebox(0,0)[r]{$#3$}}}

\def\adjust[#1]{} % For compatibility
 
\newcount \coefa
\newcount \coefb
\newcount \coefc
\newcount\tempcounta
\newcount\tempcountb
\newcount\tempcountc
\newcount\tempcountd
\newcount\xext
\newcount\yext
\newcount\xoff
\newcount\yoff
\newcount\gap%
\newcount\arrowtypea
\newcount\arrowtypeb
\newcount\arrowtypec
\newcount\arrowtyped
\newcount\arrowtypee
\newcount\height
\newcount\width
\newcount\xpos
\newcount\ypos
\newcount\run
\newcount\rise
\newcount\arrowlength
\newcount\halflength
\newcount\arrowtype
\newdimen\tempdimen
\newdimen\xlen
\newdimen\ylen
\newsavebox{\tempboxa}%
\newsavebox{\tempboxb}%
\newsavebox{\tempboxc}%
 
\newdimen\w@dth
 
\def\setw@dth#1#2{\setbox\z@\hbox{\m@th$#1$}\w@dth=\wd\z@
\setbox\@ne\hbox{\m@th$#2$}\ifnum\w@dth<\wd\@ne \w@dth=\wd\@ne \fi
\advance\w@dth by 1.2em}
 
\def\t@^#1_#2{\allowbreak\def\n@one{#1}\def\n@two{#2}\mathrel
{\setw@dth{#1}{#2}
\mathop{\hbox to \w@dth{\rightarrowfill}}\limits
\ifx\n@one\empty\else ^{\box\z@}\fi
\ifx\n@two\empty\else _{\box\@ne}\fi}}
\def\t@@^#1{\@ifnextchar_{\t@^{#1}}{\t@^{#1}_{}}}
\def\to{\@ifnextchar^{\t@@}{\t@@^{}}}
 
\def\t@left^#1_#2{\def\n@one{#1}\def\n@two{#2}\mathrel{\setw@dth{#1}{#2}
\mathop{\hbox to \w@dth{\leftarrowfill}}\limits
\ifx\n@one\empty\else ^{\box\z@}\fi
\ifx\n@two\empty\else _{\box\@ne}\fi}}
\def\t@@left^#1{\@ifnextchar_{\t@left^{#1}}{\t@left^{#1}_{}}}
\def\toleft{\@ifnextchar^{\t@@left}{\t@@left^{}}}
 
\def\two@^#1_#2{\allowbreak
\def\n@one{#1}\def\n@two{#2}\mathrel{\setw@dth{#1}{#2}
\mathop{\vcenter{\lineskip\z@\baselineskip\z@
                 \hbox to \w@dth{\rightarrowfill}%
                 \hbox to \w@dth{\rightarrowfill}}%
       }\limits
\ifx\n@one\empty\else ^{\box\z@}\fi
\ifx\n@two\empty\else _{\box\@ne}\fi}}
\def\tw@@^#1{\@ifnextchar _{\two@^{#1}}{\two@^{#1}_{}}}
\def\two{\@ifnextchar ^{\tw@@}{\tw@@^{}}}
 
\def\tofr@^#1_#2{\def\n@one{#1}\def\n@two{#2}\mathrel{\setw@dth{#1}{#2}
\mathop{\vcenter{\hbox to \w@dth{\rightarrowfill}\kern-1.7ex
                 \hbox to \w@dth{\leftarrowfill}}%
       }\limits
\ifx\n@one\empty\else ^{\box\z@}\fi
\ifx\n@two\empty\else _{\box\@ne}\fi}}
\def\t@fr@^#1{\@ifnextchar_ {\tofr@^{#1}}{\tofr@^{#1}_{}}}
\def\tofro{\@ifnextchar^ {\t@fr@}{\t@fr@^{}}}

\def\mon{\mathop{\m@th\hbox to
      14.6\P@{\lasyb\char'51\hskip-2.1\P@$\arrext$\hss
$\mathord\rightarrow$}}\limits} % width of \epi
\def\leftmono{\mathrel{\m@th\hbox to
14.6\P@{$\mathord\leftarrow$\hss$\arrext$\hskip-2.1\P@\lasyb\char'50%
}}\limits} % width of \epi
\mathchardef\arrext="0200       % amr minus for arrow extension (see \into)

\def\settypes(#1,#2,#3){\arrowtypea#1 \arrowtypeb#2 \arrowtypec#3}
\def\settoheight#1#2{\setbox\@tempboxa\hbox{#2}#1\ht\@tempboxa\relax}%
\def\settodepth#1#2{\setbox\@tempboxa\hbox{#2}#1\dp\@tempboxa\relax}%
\def\settokens`#1`#2`#3`#4`{%
     \def\tokena{#1}\def\tokenb{#2}\def\tokenc{#3}\def\tokend{#4}}
\def\setsqparms[#1`#2`#3`#4;#5`#6]{%
\arrowtypea #1
\arrowtypeb #2
\arrowtypec #3
\arrowtyped #4
\width #5
\height #6
}
\def\setpos(#1,#2){\xpos=#1 \ypos#2}

\def\settriparms[#1`#2`#3;#4]{\settripairparms[#1`#2`#3`1`1;#4]}%
 
\def\settripairparms[#1`#2`#3`#4`#5;#6]{%
\arrowtypea #1
\arrowtypeb #2
\arrowtypec #3
\arrowtyped #4
\arrowtypee #5
\width #6
\height #6
}
 
\def\resetparms{\settripairparms[1`1`1`1`1;500]\width 500}%default values%
 
\resetparms
 
\def\mvector(#1,#2)#3{%%
\put(0,0){\vector(#1,#2){#3}}%
\put(0,0){\vector(#1,#2){26}}%
}
\def\evector(#1,#2)#3{{%%
\arrowlength #3
\put(0,0){\vector(#1,#2){\arrowlength}}%
\advance \arrowlength by-30
\put(0,0){\vector(#1,#2){\arrowlength}}%
}}
 
\def\horsize#1#2{%
\settowidth{\tempdimen}{$#2$}%
#1=\tempdimen
\divide #1 by\unitlength
}
 
\def\vertsize#1#2{%
\settoheight{\tempdimen}{$#2$}%
#1=\tempdimen
\settodepth{\tempdimen}{$#2$}%
\advance #1 by\tempdimen
\divide #1 by\unitlength
}
 
\def\putvector(#1,#2)(#3,#4)#5#6{{%
\ifnum3<\arrowtype
\putdashvector(#1,#2)(#3,#4)#5\arrowtype
\else
\ifnum\arrowtype<-3
\putdashvector(#1,#2)(#3,#4)#5\arrowtype
\else
\xpos=#1
\ypos=#2
\run=#3
\rise=#4
\arrowlength=#5
\ifnum \arrowtype<0
    \ifnum \run=0
        \advance \ypos by-\arrowlength
    \else
        \tempcounta \arrowlength
        \multiply \tempcounta by\rise
        \divide \tempcounta by\run
        \ifnum\run>0
            \advance \xpos by\arrowlength
            \advance \ypos by\tempcounta
        \else
            \advance \xpos by-\arrowlength
            \advance \ypos by-\tempcounta
        \fi
    \fi
    \multiply \arrowtype by-1
    \multiply \rise by-1
    \multiply \run by-1
\fi
\ifcase \arrowtype
\or \put(\xpos,\ypos){\vector(\run,\rise){\arrowlength}}%
\or \put(\xpos,\ypos){\mvector(\run,\rise)\arrowlength}%
\or \put(\xpos,\ypos){\evector(\run,\rise){\arrowlength}}%
\fi\fi\fi
}}
 
\def\putsplitvector(#1,#2)#3#4{%%
\xpos #1
\ypos #2
\arrowtype #4
\halflength #3
\arrowlength #3
\gap 140
\advance \halflength by-\gap
\divide \halflength by2
\ifnum\arrowtype>0
   \ifcase \arrowtype
   \or \put(\xpos,\ypos){\line(0,-1){\halflength}}%
       \advance\ypos by-\halflength
       \advance\ypos by-\gap
       \put(\xpos,\ypos){\vector(0,-1){\halflength}}%
   \or \put(\xpos,\ypos){\line(0,-1)\halflength}%
       \put(\xpos,\ypos){\vector(0,-1)3}%
       \advance\ypos by-\halflength
       \advance\ypos by-\gap
       \put(\xpos,\ypos){\vector(0,-1){\halflength}}%
   \or \put(\xpos,\ypos){\line(0,-1)\halflength}%
       \advance\ypos by-\halflength
       \advance\ypos by-\gap
       \put(\xpos,\ypos){\evector(0,-1){\halflength}}%
   \fi
\else \arrowtype=-\arrowtype
   \ifcase\arrowtype
   \or \advance \ypos by-\arrowlength
       \put(\xpos,\ypos){\line(0,1){\halflength}}%
       \advance\ypos by\halflength
       \advance\ypos by\gap
       \put(\xpos,\ypos){\vector(0,1){\halflength}}%
   \or \advance \ypos by-\arrowlength
       \put(\xpos,\ypos){\line(0,1)\halflength}%
       \put(\xpos,\ypos){\vector(0,1)3}%
       \advance\ypos by\halflength
       \advance\ypos by\gap
       \put(\xpos,\ypos){\vector(0,1){\halflength}}%
   \or \advance \ypos by-\arrowlength
       \put(\xpos,\ypos){\line(0,1)\halflength}%
       \advance\ypos by\halflength
       \advance\ypos by\gap
       \put(\xpos,\ypos){\evector(0,1){\halflength}}%
   \fi
\fi
}
 
\def\putmorphism(#1)(#2,#3)[#4`#5`#6]#7#8#9{{%
\run #2
\rise #3
\ifnum\rise=0
  \puthmorphism(#1)[#4`#5`#6]{#7}{#8}#9%
\else\ifnum\run=0
  \putvmorphism(#1)[#4`#5`#6]{#7}{#8}#9%
\else
\setpos(#1)%
\arrowlength #7
\arrowtype #8
\ifnum\run=0
\else\ifnum\rise=0
\else
\ifnum\run>0
    \coefa=1
\else
   \coefa=-1
\fi
\ifnum\arrowtype>0
   \coefb=0
   \coefc=-1
\else
   \coefb=\coefa
   \coefc=1
   \arrowtype=-\arrowtype
\fi
\width=2
\multiply \width by\run
\divide \width by\rise
\ifnum \width<0  \width=-\width\fi
\advance\width by60
\if l#9 \width=-\width\fi
\putbox(\xpos,\ypos){#4}%            %node 1
{\multiply \coefa by\arrowlength%      %node 2
\advance\xpos by\coefa
\multiply \coefa by\rise
\divide \coefa by\run
\advance \ypos by\coefa
\putbox(\xpos,\ypos){#5} }%
{\multiply \coefa by\arrowlength%      %label
\divide \coefa by2
\advance \xpos by\coefa
\advance \xpos by\width
\multiply \coefa by\rise
\divide \coefa by\run
\advance \ypos by\coefa
\if l#9%
   \putrbox(\xpos,\ypos){#6}%
\else\if r#9%
   \putlbox(\xpos,\ypos){#6}%
\fi\fi }%
{\multiply \rise by-\coefc%             %arrow
\multiply \run by-\coefc
\multiply \coefb by\arrowlength
\advance \xpos by\coefb
\multiply \coefb by\rise
\divide \coefb by\run
\advance \ypos by\coefb
\multiply \coefc by70
\advance \ypos by\coefc
\multiply \coefc by\run
\divide \coefc by\rise
\advance \xpos by\coefc
\multiply \coefa by140
\multiply \coefa by\run
\divide \coefa by\rise
\advance \arrowlength by\coefa
\ifcase\arrowtype
\or \put(\xpos,\ypos){\vector(\run,\rise){\arrowlength}}%
\or \put(\xpos,\ypos){\mvector(\run,\rise){\arrowlength}}%
\or \put(\xpos,\ypos){\evector(\run,\rise){\arrowlength}}%
\fi}\fi\fi\fi\fi}}

\newcount\numbdashes \newcount\lengthdash \newcount\increment
 
\def\howmanydashes{% Actually returns both number and length
\numbdashes=\arrowlength \lengthdash=40
\divide\numbdashes by \lengthdash
\lengthdash=\arrowlength
\divide\lengthdash by \numbdashes
\increment=\lengthdash
\multiply\lengthdash by 3
\divide\lengthdash by 5
}
 
\def\putdashvector(#1)(#2,#3)#4#5{%
\ifnum#3=0 \putdashhvector(#1){#4}#5
\else
\ifnum#2=0
\putdashvvector(#1){#4}#5\fi\fi}
 
\def\putdashhvector(#1,#2)#3#4{{%
\arrowlength=#3 \howmanydashes
\multiput(#1,#2)(\increment,0){\numbdashes}%
{\vrule height .4pt width \lengthdash\unitlength}
\arrowtype=#4 \xpos=#1
\ifnum\arrowtype<0 \advance\arrowtype by 7 \fi
\ifcase\arrowtype
\or \advance\xpos by 10
    \put(\xpos,#2){\vector(-1,0){\lengthdash}}
    \advance\xpos by 40
    \put(\xpos,#2){\vector(-1,0){\lengthdash}}
\or \advance \xpos by 10
    \put(\xpos,#2){\vector(-1,0){\lengthdash}}
    \advance\xpos by  \arrowlength
    \advance\xpos by  -50
    \put(\xpos,#2){\vector(-1,0){\lengthdash}}
\or \advance\xpos by 10
    \put(\xpos,#2){\vector(-1,0){\lengthdash}}
\or \advance\xpos by \arrowlength
    \advance\xpos by -\lengthdash
    \put(\xpos,#2){\vector(1,0){\lengthdash}}
\or {\advance\xpos by 10
    \put(\xpos,#2){\vector(1,0){\lengthdash}}}
    \advance\xpos by \arrowlength
    \advance\xpos by -\lengthdash
    \put(\xpos,#2){\vector(1,0){\lengthdash}}
\or \advance\xpos by \arrowlength
    \advance\xpos by -\lengthdash
    \put(\xpos,#2){\vector(1,0){\lengthdash}}
    \advance\xpos by -40
    \put(\xpos,#2){\vector(1,0){\lengthdash}}
   \fi
}}
 
\def\putdashvvector(#1,#2)#3#4{{%
\arrowlength=#3 \howmanydashes
\ypos=#2 \advance\ypos by -\arrowlength
\multiput(#1,#2)(0,\increment){\numbdashes}%
    {\vrule width .4pt height \lengthdash\unitlength}
\arrowtype=#4 \ypos=#2
\ifnum\arrowtype<0 \advance\arrowtype by 7 \fi
\ifcase\arrowtype
\or \advance\ypos by \arrowlength \advance\ypos by -40
    \put(#1,\ypos){\vector(0,1){\lengthdash}}
    \advance\ypos by -40
    \put(#1,\ypos){\vector(0,1){\lengthdash}}
\or \advance\ypos by 10
    \put(#1,\ypos){\vector(0,1){\lengthdash}}
    \advance\ypos by \arrowlength \advance\ypos by -40
    \put(#1,\ypos){\vector(0,1){\lengthdash}}
\or \advance\ypos by \arrowlength \advance\ypos by -40
    \put(#1,\ypos){\vector(0,1){\lengthdash}}
\or \advance\ypos by 10
    \put(#1,\ypos){\vector(0,-1){\lengthdash}}
\or \advance\ypos by 10
    \put(#1,\ypos){\vector(0,-1){\lengthdash}}
    \advance\ypos by \arrowlength \advance\ypos by -40
    \put(#1,\ypos){\vector(0,-1){\lengthdash}}
\or \advance\ypos by 10
    \put(#1,\ypos){\vector(0,-1){\lengthdash}}
    \advance\ypos by 40
    \put(#1,\ypos){\vector(0,-1){\lengthdash}}
\fi
}}
 
\def\puthmorphism(#1,#2)[#3`#4`#5]#6#7#8{{%
\xpos #1
\ypos #2
\width #6
\arrowlength #6
\arrowtype=#7
\putbox(\xpos,\ypos){#3\vphantom{#4}}%
{\advance \xpos by\arrowlength
\putbox(\xpos,\ypos){\vphantom{#3}#4}}%
\horsize{\tempcounta}{#3}%
\horsize{\tempcountb}{#4}%
\divide \tempcounta by2
\divide \tempcountb by2
\advance \tempcounta by30
\advance \tempcountb by30
\advance \xpos by\tempcounta
\advance \arrowlength by-\tempcounta
\advance \arrowlength by-\tempcountb
\putvector(\xpos,\ypos)(1,0)\arrowlength\arrowtype
\divide \arrowlength by2
\advance \xpos by\arrowlength
\vertsize{\tempcounta}{#5}%
\divide\tempcounta by2
\advance \tempcounta by20
\if a#8 %
   \advance \ypos by\tempcounta
   \putbox(\xpos,\ypos){#5}%
\else
   \advance \ypos by-\tempcounta
   \putbox(\xpos,\ypos){#5}%
\fi}}
 
\def\putvmorphism(#1,#2)[#3`#4`#5]#6#7#8{{%
\xpos #1
\ypos #2
\arrowlength #6
\arrowtype #7
\settowidth{\xlen}{$#5$}%
\putbox(\xpos,\ypos){#3}%
{\advance \ypos by-\arrowlength
\putbox(\xpos,\ypos){#4}}%
{\advance\arrowlength by-140
\advance \ypos by-70
\ifdim\xlen>0pt
   \if m#8%
      \putsplitvector(\xpos,\ypos)\arrowlength\arrowtype
   \else
   \putvector(\xpos,\ypos)(0,-1)\arrowlength\arrowtype
   \fi
\else
   \putvector(\xpos,\ypos)(0,-1)\arrowlength\arrowtype
\fi}%
\ifdim\xlen>0pt
   \divide \arrowlength by2
   \advance\ypos by-\arrowlength
   \if l#8%
      \advance \xpos by-40
      \putrbox(\xpos,\ypos){#5}%
   \else\if r#8%
      \advance \xpos by40
      \putlbox(\xpos,\ypos){#5}%
   \else
      \putbox(\xpos,\ypos){#5}%
   \fi\fi
\fi
}}
 
\def\putsquarep<#1>(#2)[#3;#4`#5`#6`#7]{{%
\setsqparms[#1]%
\setpos(#2)%
\settokens`#3`%
\puthmorphism(\xpos,\ypos)[\tokenc`\tokend`{#7}]{\width}{\arrowtyped}b%
\advance\ypos by \height
\puthmorphism(\xpos,\ypos)[\tokena`\tokenb`{#4}]{\width}{\arrowtypea}a%
\putvmorphism(\xpos,\ypos)[``{#5}]{\height}{\arrowtypeb}l%
\advance\xpos by \width
\putvmorphism(\xpos,\ypos)[``{#6}]{\height}{\arrowtypec}r%
}}
 
\def\putsquare{\@ifnextchar <{\putsquarep}{\putsquarep%
   <\arrowtypea`\arrowtypeb`\arrowtypec`\arrowtyped;\width`\height>}}
\def\square{\@ifnextchar< {\squarep}{\squarep
   <\arrowtypea`\arrowtypeb`\arrowtypec`\arrowtyped;\width`\height>}}
\def\squarep<#1>[#2`#3`#4`#5;#6`#7`#8`#9]{{%       %     #2------>#3
\setsqparms[#1]%                                   %      |       |
\diagram%                                          %      |       |
\putsquarep<\arrowtypea`\arrowtypeb`\arrowtypec`%  %    #7|       |#8
\arrowtyped;\width`\height>%                       %      |       |
(0,0)[#2`#3`#4`{#5};#6`#7`#8`{#9}]%                %      |       |
\enddiagram%                                       %      v       v
}}                                                 %     #4------>#5
\def\putptrianglep<#1>(#2,#3)[#4`#5`#6;#7`#8`#9]{{%
\settriparms[#1]%
\xpos=#2 \ypos=#3
\advance\ypos by \height
\puthmorphism(\xpos,\ypos)[#4`#5`{#7}]{\height}{\arrowtypea}a%
\putvmorphism(\xpos,\ypos)[`#6`{#8}]{\height}{\arrowtypeb}l%
\advance\xpos by\height
\putmorphism(\xpos,\ypos)(-1,-1)[``{#9}]{\height}{\arrowtypec}r%
}}
 
\def\putptriangle{\@ifnextchar <{\putptrianglep}{\putptrianglep
   <\arrowtypea`\arrowtypeb`\arrowtypec;\height>}}
\def\ptriangle{\@ifnextchar <{\ptrianglep}{\ptrianglep
   <\arrowtypea`\arrowtypeb`\arrowtypec;\height>}}
\def\ptrianglep<#1>[#2`#3`#4;#5`#6`#7]{{%%    %      #2----->#3
\settriparms[#1]%                             %      |      /
\diagram%                                     %      |     /
\putptrianglep<\arrowtypea`\arrowtypeb`%      %    #6|    /#7
\arrowtypec;\height>%                         %      |   /
(0,0)[#2`#3`#4;#5`#6`{#7}]%                   %      |  /
\enddiagram%%                                 %      v v
}}                                            %      #4
 
\def\putqtrianglep<#1>(#2,#3)[#4`#5`#6;#7`#8`#9]{{%
\settriparms[#1]%
\xpos=#2 \ypos=#3
\advance\ypos by\height
\puthmorphism(\xpos,\ypos)[#4`#5`{#7}]{\height}{\arrowtypea}a%
\putmorphism(\xpos,\ypos)(1,-1)[``{#8}]{\height}{\arrowtypeb}l%
\advance\xpos by\height
\putvmorphism(\xpos,\ypos)[`#6`{#9}]{\height}{\arrowtypec}r%
}}
 
\def\putqtriangle{\@ifnextchar <{\putqtrianglep}{\putqtrianglep
   <\arrowtypea`\arrowtypeb`\arrowtypec;\height>}}
\def\qtriangle{\@ifnextchar <{\qtrianglep}{\qtrianglep
   <\arrowtypea`\arrowtypeb`\arrowtypec;\height>}}
\def\qtrianglep<#1>[#2`#3`#4;#5`#6`#7]{{%%    %        #2----->#3
\settriparms[#1]%                             %         \      |
\width=\height                                %          \     |
\diagram%                                     %         #6\    |#7
\putqtrianglep<\arrowtypea`\arrowtypeb`%      %            \   |
\arrowtypec;\height>%                         %             \  |
(0,0)[#2`#3`#4;#5`#6`{#7}]%                   %              v v
\enddiagram%%                                 %               #4
}}
 
\def\putdtrianglep<#1>(#2,#3)[#4`#5`#6;#7`#8`#9]{{%
\settriparms[#1]%
\xpos=#2 \ypos=#3
\puthmorphism(\xpos,\ypos)[#5`#6`{#9}]{\height}{\arrowtypec}b%
\advance\xpos by \height \advance\ypos by\height
\putmorphism(\xpos,\ypos)(-1,-1)[``{#7}]{\height}{\arrowtypea}l%
\putvmorphism(\xpos,\ypos)[#4``{#8}]{\height}{\arrowtypeb}r%
}}
 
\def\putdtriangle{\@ifnextchar <{\putdtrianglep}{\putdtrianglep
   <\arrowtypea`\arrowtypeb`\arrowtypec;\height>}}
\def\dtriangle{\@ifnextchar <{\dtrianglep}{\dtrianglep
   <\arrowtypea`\arrowtypeb`\arrowtypec;\height>}}
\def\dtrianglep<#1>[#2`#3`#4;#5`#6`#7]{{%%    %                  / |
\settriparms[#1]%                             %                 /  |
\width=\height                                %              #5/   |#6
\diagram%                                     %               /    |
\putdtrianglep<\arrowtypea`\arrowtypeb`%      %              /     |
\arrowtypec;\height>%                         %             v      v
(0,0)[#2`#3`#4;#5`#6`{#7}]%                   %            #3----->#4
\enddiagram%%                                 %                #7
}}
 
\def\putbtrianglep<#1>(#2,#3)[#4`#5`#6;#7`#8`#9]{{%
\settriparms[#1]%
\xpos=#2 \ypos=#3
\puthmorphism(\xpos,\ypos)[#5`#6`{#9}]{\height}{\arrowtypec}b%
\advance\ypos by\height
\putmorphism(\xpos,\ypos)(1,-1)[``{#8}]{\height}{\arrowtypeb}r%
\putvmorphism(\xpos,\ypos)[#4``{#7}]{\height}{\arrowtypea}l%
}}
 
\def\putbtriangle{\@ifnextchar <{\putbtrianglep}{\putbtrianglep
   <\arrowtypea`\arrowtypeb`\arrowtypec;\height>}}
\def\btriangle{\@ifnextchar <{\btrianglep}{\btrianglep
   <\arrowtypea`\arrowtypeb`\arrowtypec;\height>}}
\def\btrianglep<#1>[#2`#3`#4;#5`#6`#7]{{%%   %              | \
\settriparms[#1]%                            %              |  \
\width=\height                               %            #5|   \#6
\diagram%                                    %              |    \
\putbtrianglep<\arrowtypea`\arrowtypeb`%     %              |     \
\arrowtypec;\height>%                        %              v      v
(0,0)[#2`#3`#4;#5`#6`{#7}]%                  %              #3----->#4
\enddiagram%%                                %                 #7
}}
 
\def\putAtrianglep<#1>(#2,#3)[#4`#5`#6;#7`#8`#9]{{%
\settriparms[#1]%
\xpos=#2 \ypos=#3
{\multiply \height by2
\puthmorphism(\xpos,\ypos)[#5`#6`{#9}]{\height}{\arrowtypec}b}%
\advance\xpos by\height \advance\ypos by\height
\putmorphism(\xpos,\ypos)(-1,-1)[#4``{#7}]{\height}{\arrowtypea}l%
\putmorphism(\xpos,\ypos)(1,-1)[``{#8}]{\height}{\arrowtypeb}r%
}}
 
\def\putAtriangle{\@ifnextchar <{\putAtrianglep}{\putAtrianglep
   <\arrowtypea`\arrowtypeb`\arrowtypec;\height>}}
\def\Atriangle{\@ifnextchar <{\Atrianglep}{\Atrianglep
   <\arrowtypea`\arrowtypeb`\arrowtypec;\height>}}
\def\Atrianglep<#1>[#2`#3`#4;#5`#6`#7]{{%%         %         /   \
\settriparms[#1]%                                  %        /     \
\width=\height                                     %     #5/       \#6
\diagram%                                          %      /         \
\putAtrianglep<\arrowtypea`\arrowtypeb`%           %     /           \
\arrowtypec;\height>%                              %    v             v
(0,0)[#2`#3`#4;#5`#6`{#7}]%                        %   #3------------>#4
\enddiagram%%                                      %          #7
}}
 
\def\putAtrianglepairp<#1>(#2)[#3;#4`#5`#6`#7`#8]{{%
\settripairparms[#1]%
\setpos(#2)%
\settokens`#3`%
\puthmorphism(\xpos,\ypos)[\tokenb`\tokenc`{#7}]{\height}{\arrowtyped}b%
\advance\xpos by\height
\puthmorphism(\xpos,\ypos)[\phantom{\tokenc}`\tokend`{#8}]%
{\height}{\arrowtypee}b%
\advance\ypos by\height
\putmorphism(\xpos,\ypos)(-1,-1)[\tokena``{#4}]{\height}{\arrowtypea}l%
\putvmorphism(\xpos,\ypos)[``{#5}]{\height}{\arrowtypeb}m%
\putmorphism(\xpos,\ypos)(1,-1)[``{#6}]{\height}{\arrowtypec}r%
}}
 
\def\putAtrianglepair{\@ifnextchar <{\putAtrianglepairp}{\putAtrianglepairp%
   <\arrowtypea`\arrowtypeb`\arrowtypec`\arrowtyped`\arrowtypee;\height>}}
\def\Atrianglepair{\@ifnextchar <{\Atrianglepairp}{\Atrianglepairp%
   <\arrowtypea`\arrowtypeb`\arrowtypec`\arrowtyped`\arrowtypee;\height>}}
 
\def\Atrianglepairp<#1>[#2;#3`#4`#5`#6`#7]{{%           %  #2a
\settripairparms[#1]%                         %           / | \
\settokens`#2`%                               %          /  |  \
\width=\height                                %       #3/  #4   \#5
\diagram%                                     %        /    |    \
\putAtrianglepairp                            %       /     |     \
<\arrowtypea`\arrowtypeb`\arrowtypec`%        %      v      v      v
\arrowtyped`\arrowtypee;\height>%             %     #2b---->#2c---->#2d
(0,0)[{#2};#3`#4`#5`#6`{#7}]%                 %         #6     #7
\enddiagram%%
}}
 
\def\putVtrianglep<#1>(#2,#3)[#4`#5`#6;#7`#8`#9]{{%
\settriparms[#1]%
\xpos=#2 \ypos=#3
\advance\ypos by\height
{\multiply\height by2
\puthmorphism(\xpos,\ypos)[#4`#5`{#7}]{\height}{\arrowtypea}a}%
\putmorphism(\xpos,\ypos)(1,-1)[`#6`{#8}]{\height}{\arrowtypeb}l%
\advance\xpos by\height
\advance\xpos by\height
\putmorphism(\xpos,\ypos)(-1,-1)[``{#9}]{\height}{\arrowtypec}r%
}}
 
\def\putVtriangle{\@ifnextchar <{\putVtrianglep}{\putVtrianglep
   <\arrowtypea`\arrowtypeb`\arrowtypec;\height>}}
\def\Vtriangle{\@ifnextchar <{\Vtrianglep}{\Vtrianglep
   <\arrowtypea`\arrowtypeb`\arrowtypec;\height>}}
\def\Vtrianglep<#1>[#2`#3`#4;#5`#6`#7]{{%%     %        #2------------->#3
\settriparms[#1]%                              %         \             /
\width=\height                                 %          \           /
\diagram%                                      %         #6\         /#7
\putVtrianglep<\arrowtypea`\arrowtypeb`%       %            \       /
\arrowtypec;\height>%                          %             \     /
(0,0)[#2`#3`#4;#5`#6`{#7}]%                    %              v   v
\enddiagram%%                                  %               #4
}}
 
\def\putVtrianglepairp<#1>(#2)[#3;#4`#5`#6`#7`#8]{{
\settripairparms[#1]%
\setpos(#2)%
\settokens`#3`%
\advance\ypos by\height
\putmorphism(\xpos,\ypos)(1,-1)[`\tokend`{#6}]{\height}{\arrowtypec}l%
\puthmorphism(\xpos,\ypos)[\tokena`\tokenb`{#4}]{\height}{\arrowtypea}a%
\advance\xpos by\height
\puthmorphism(\xpos,\ypos)[\phantom{\tokenb}`\tokenc`{#5}]%
{\height}{\arrowtypeb}a%
\putvmorphism(\xpos,\ypos)[``{#7}]{\height}{\arrowtyped}m%
\advance\xpos by\height
\putmorphism(\xpos,\ypos)(-1,-1)[``{#8}]{\height}{\arrowtypee}r%
}}
 
\def\putVtrianglepair{\@ifnextchar <{\putVtrianglepairp}{\putVtrianglepairp%
    <\arrowtypea`\arrowtypeb`\arrowtypec`\arrowtyped`\arrowtypee;\height>}}
\def\Vtrianglepair{\@ifnextchar <{\Vtrianglepairp}{\Vtrianglepairp%
    <\arrowtypea`\arrowtypeb`\arrowtypec`\arrowtyped`\arrowtypee;\height>}}
\def\Vtrianglepairp<#1>[#2;#3`#4`#5`#6`#7]{{%  %  #2a---->#2b---->#2c
\settripairparms[#1]%                          %   \      |      /
\settokens`#2`%                                %    \     |     /
\diagram%                                      %   #5\   #6    /#7
\putVtrianglepairp                             %      \   |   /
<\arrowtypea`\arrowtypeb`\arrowtypec`%         %       \  |  /
\arrowtyped`\arrowtypee;\height>%              %        v v v
(0,0)[{#2};#3`#4`#5`#6`{#7}]%                  %         #2d
\enddiagram%%
}}

\def\putCtrianglep<#1>(#2,#3)[#4`#5`#6;#7`#8`#9]{{%
\settriparms[#1]%
\xpos=#2 \ypos=#3
\advance\ypos by\height
\putmorphism(\xpos,\ypos)(1,-1)[``{#9}]{\height}{\arrowtypec}l%
\advance\xpos by\height
\advance\ypos by\height
\putmorphism(\xpos,\ypos)(-1,-1)[#4`#5`{#7}]{\height}{\arrowtypea}l%
{\multiply\height by 2
\putvmorphism(\xpos,\ypos)[`#6`{#8}]{\height}{\arrowtypeb}r}%
}}
 
\def\putCtriangle{\@ifnextchar <{\putCtrianglep}{\putCtrianglep
    <\arrowtypea`\arrowtypeb`\arrowtypec;\height>}}
\def\Ctriangle{\@ifnextchar <{\Ctrianglep}{\Ctrianglep
    <\arrowtypea`\arrowtypeb`\arrowtypec;\height>}}
\def\Ctrianglep<#1>[#2`#3`#4;#5`#6`#7]{{%%   %                / |
\settriparms[#1]%                            %             #5/  |
\width=\height                               %              /   |
\diagram%                                    %             v    |
\putCtrianglep<\arrowtypea`\arrowtypeb`%     %           #3     |#6
\arrowtypec;\height>%                        %             \    |
(0,0)[#2`#3`#4;#5`#6`{#7}]%                  %            #7\   |
\enddiagram%%                                %               \  |
}}                                           %                v v
\def\putDtrianglep<#1>(#2,#3)[#4`#5`#6;#7`#8`#9]{{%
\settriparms[#1]%
\xpos=#2 \ypos=#3
\advance\xpos by\height \advance\ypos by\height
\putmorphism(\xpos,\ypos)(-1,-1)[``{#9}]{\height}{\arrowtypec}r%
\advance\xpos by-\height \advance\ypos by\height
\putmorphism(\xpos,\ypos)(1,-1)[`#5`{#8}]{\height}{\arrowtypeb}r%
{\multiply\height by 2
\putvmorphism(\xpos,\ypos)[#4`#6`{#7}]{\height}{\arrowtypea}l}%
}}
 
\def\putDtriangle{\@ifnextchar <{\putDtrianglep}{\putDtrianglep
    <\arrowtypea`\arrowtypeb`\arrowtypec;\height>}}
\def\Dtriangle{\@ifnextchar <{\Dtrianglep}{\Dtrianglep
   <\arrowtypea`\arrowtypeb`\arrowtypec;\height>}}
\def\Dtrianglep<#1>[#2`#3`#4;#5`#6`#7]{{%%  %          | \
\settriparms[#1]%                           %          |  \#6
\width=\height                              %          |   \
\diagram%                                   %          |    v
\putDtrianglep<\arrowtypea`\arrowtypeb`%    %        #5|    #3
\arrowtypec;\height>%                       %          |    /
(0,0)[#2`#3`#4;#5`#6`{#7}]%                 %          |   /#7
\enddiagram%%                               %          |  /
}}                                          %          v v
\def\setrecparms[#1`#2]{\width=#1 \height=#2}%
\def\recursep<#1`#2>[#3;#4`#5`#6`#7`#8]{{\m@th
\width=#1 \height=#2
\settokens`#3`
\settowidth{\tempdimen}{$\tokena$}
\ifdim\tempdimen=0pt
  \savebox{\tempboxa}{\hbox{$\tokenb$}}%
  \savebox{\tempboxb}{\hbox{$\tokend$}}%
  \savebox{\tempboxc}{\hbox{$#6$}}%
\else
  \savebox{\tempboxa}{\hbox{$\hbox{$\tokena$}\times\hbox{$\tokenb$}$}}%
  \savebox{\tempboxb}{\hbox{$\hbox{$\tokena$}\times\hbox{$\tokend$}$}}%
  \savebox{\tempboxc}{\hbox{$\hbox{$\tokena$}\times\hbox{$#6$}$}}%
\fi
\ypos=\height
\divide\ypos by 2
\xpos=\ypos
\advance\xpos by \width
\bfig
\putCtrianglep<-1`1`1;\ypos>(0,0)[`\tokenc`;#5`#6`{#7}]%
\puthmorphism(\ypos,0)[\tokend`\usebox{\tempboxb}`{#8}]{\width}{-1}b%
\puthmorphism(\ypos,\height)[\tokenb`\usebox{\tempboxa}`{#4}]{\width}{-1}a%
\advance\ypos by \width
\putvmorphism(\ypos,\height)[``\usebox{\tempboxc}]{\height}1r%
\efig
}}
 
\def\recurse{\@ifnextchar <{\recursep}{\recursep<\width`\height>}}
 
\def\puttwohmorphisms(#1,#2)[#3`#4;#5`#6]#7#8#9{{%
\puthmorphism(#1,#2)[#3`#4`]{#7}0a
\ypos=#2
\advance\ypos by 20
\puthmorphism(#1,\ypos)[\phantom{#3}`\phantom{#4}`#5]{#7}{#8}a
\advance\ypos by -40
\puthmorphism(#1,\ypos)[\phantom{#3}`\phantom{#4}`#6]{#7}{#9}b
}}
 
\def\puttwovmorphisms(#1,#2)[#3`#4;#5`#6]#7#8#9{{%
\putvmorphism(#1,#2)[#3`#4`]{#7}0a
\xpos=#1
\advance\xpos by -20
\putvmorphism(\xpos,#2)[\phantom{#3}`\phantom{#4}`#5]{#7}{#8}l
\advance\xpos by 40
\putvmorphism(\xpos,#2)[\phantom{#3}`\phantom{#4}`#6]{#7}{#9}r
}}
 
\def\puthcoequalizer(#1)[#2`#3`#4;#5`#6`#7]#8#9{{%
\setpos(#1)%
\puttwohmorphisms(\xpos,\ypos)[#2`#3;#5`#6]{#8}11%
\advance\xpos by #8
\puthmorphism(\xpos,\ypos)[\phantom{#3}`#4`#7]{#8}1{#9}
}}
 
\def\putvcoequalizer(#1)[#2`#3`#4;#5`#6`#7]#8#9{{%
\setpos(#1)%
\puttwovmorphisms(\xpos,\ypos)[#2`#3;#5`#6]{#8}11%
\advance\ypos by -#8
\putvmorphism(\xpos,\ypos)[\phantom{#3}`#4`#7]{#8}1{#9}
}}
 
\def\putthreehmorphisms(#1)[#2`#3;#4`#5`#6]#7(#8)#9{{%
\setpos(#1) \settypes(#8)
\if a#9 %
     \vertsize{\tempcounta}{#5}%
     \vertsize{\tempcountb}{#6}%
     \ifnum \tempcounta<\tempcountb \tempcounta=\tempcountb \fi
\else
     \vertsize{\tempcounta}{#4}%
     \vertsize{\tempcountb}{#5}%
     \ifnum \tempcounta<\tempcountb \tempcounta=\tempcountb \fi
\fi
\advance \tempcounta by 60
\puthmorphism(\xpos,\ypos)[#2`#3`#5]{#7}{\arrowtypeb}{#9}
\advance\ypos by \tempcounta
\puthmorphism(\xpos,\ypos)[\phantom{#2}`\phantom{#3}`#4]{#7}{\arrowtypea}{#9}
\advance\ypos by -\tempcounta \advance\ypos by -\tempcounta
\puthmorphism(\xpos,\ypos)[\phantom{#2}`\phantom{#3}`#6]{#7}{\arrowtypec}{#9}
}}
 
\def\setarrowtoks[#1`#2`#3`#4`#5`#6]{%
\def\toka{#1}
\def\tokb{#2}
\def\tokc{#3}
\def\tokd{#4}
\def\toke{#5}
\def\tokf{#6}
}
\def\hex{\@ifnextchar <{\hexp}{\hexp<1000`400>}}
\def\hexp<#1`#2>[#3`#4`#5`#6`#7`#8;#9]{%
\setarrowtoks[#9]
\yext=#2 \advance \yext by #2
\xext=#1 \advance\xext by \yext
\bfig
\putCtriangle<-1`0`1;#2>(0,0)[`#5`;\tokb``\tokd]
\xext=#1 \yext=#2 \advance \yext by #2
\putsquare<1`0`0`1;\xext`\yext>(#2,0)[#3`#4`#7`#8;\toka```\tokf]
\advance \xext by #2
\putDtriangle<0`1`-1;#2>(\xext,0)[`#6`;`\tokc`\toke]
\efig
}
\makeatother

\newtheorem{thm}{Theorem}
\newtheorem{prop}[thm]{Proposition}
\newtheorem{cor}[thm]{Corollary}
\newtheorem{lem}[thm]{Lemma}

\theoremstyle{remark}
\newtheorem{rem}[thm]{Remark}

\newcommand{\C}{\mathbb{C}}

\newcommand{\ga}{\alpha}

\newcommand{\gc}{\gamma}

\newcommand{\inv}{^{-1}}

\newcommand{\Q}{\mathbb{Q}}

\newcommand{\Sym}{\operatorname{Sym}}

\newcommand{\tri}{\bigtriangleup}

\begin{document}
\title[Localization formulas in equivariant cohomology]
{On the localization formula in equivariant
cohomology} 

\author{Andr\'es Pedroza}
\address{Facultad de Ciencias\\
         Universidad de Colima\\
	 Bernal D\'{\i}az del Castillo 340\\
         Colima, Colima, Mexico 28045}
\email{andres\_pedroza@ucol.mx}

\author{Loring W. Tu}
\address{Department of Mathematics\\
         Tufts University\\
         Medford, MA, USA 02155}
\email{loring.tu@tufts.edu}
\thanks{The first author was supported in part by FRABA-Universidad
de Colima grant.
The second author acknowledges the hospitality and support of
the Institut Henri Poincar\'e
the Institut de Math\'ematiques de Jussieu, Paris.}

\date{\today}
\begin{abstract}
We give a generalization of the Atiyah-Bott-Berline-Vergne localization theorem for the
equivariant cohomology of a torus action. We
replace the manifold having a torus action by an equivariant
map of manifolds having a compact connected Lie group action.  This provides a systematic method for calculating the
Gysin homomorphism in ordinary cohomology of an equivariant
map.  As an example, we recover a formula of
Akyildiz-Carrell for the Gysin homomorphism of flag
manifolds.

\end{abstract}
\maketitle

Suppose $M$ is a compact oriented manifold on which a torus $T$
acts.  The Atiyah-Bott-Berline-Vergne localization formula
calculates the integral of an equivariant cohomology class on $M$
in terms of an integral over the fixed point set $M^T$. This
formula has found many applications, for example, in analysis,
topology, symplectic geometry, and algebraic geometry (see
\cite{atiyah-bott}, \cite{duistermaat}, \cite{ellingsrud-stromme},
\cite{tu-gysin}). Similar, but not entirely analogous, formulas
exist in $K$-theory (\cite{atiyah-segal}), cobordism theory
(\cite{quillen}), and algebraic geometry (\cite{edidin-graham}).

Taking cues from the work of Atiyah and Segal in $K$-theory
\cite{atiyah-segal}, we state
and prove a localization formula for a compact connected Lie
group action in terms of the
fixed point set of a conjugacy class in the group.
As an application, the formula can be used to calculate
the Gysin homomorphism in ordinary cohomology of an equivariant
map.
For a compact connected Lie group $G$ with maximal torus
$T$ and a closed subgroup $H$ containing $T$, we work out
as an example the Gysin homomorphism of the canonical
projection $f: G/T \to G/H$, a formula first obtained by
Akyildiz and Carrell \cite{akyildiz-carrell}.

The application to the Gysin map in this article complements
that of \cite{tu-gysin}.
The previous article \cite{tu-gysin} shows how to use the
ABBV localization formula to calculate the Gysin map of a
fiber bundle.  This article shows how to use the relative
localization formula to calculate the Gysin map of an
equivariant map.

We thank Michel Brion for many helpful discussions.

\section{Borel-type localization formula for a conjugacy class}
\setcounter{thm}{0}

Suppose a compact connected Lie group $G$ acts on a manifold
$M$.  For $g \in G$, define $M^g$ to be the fixed point set
of $g$:
\[
M^g = \{ x \in M \mid g\cdot x = x \}.
\]
The set $M^g$ is not $G$-invariant. The $G$-invariant subset it
generates is
$$
\cup_{h\in G} h\cdot (M^g)= \cup_{h\in G} M^{hgh^{-1}}=
\cup_{k\in C(g)} M^k
$$
where $C(g)$ is the conjugacy class of $g$. This suggests that
for every conjugacy class $C$ in $G$, we consider the set $M^C$
of elements of $M$ that are fixed by at least one element
of the conjugacy class $C$:
\[
M^C = \cup_{k\in C} M^k.
\]
Then $M^C$ is a closed $G$-subset of $M$ (\cite{atiyah-segal},
footnote 1, p.~532); however it may not be always smooth. From
now on we  make the assumption that $M^C$ is smooth.

Suppose $C = C(g)$ is the conjugacy class of an element $g$ in $G$.  Let $T$ be a maximal torus of $T$ containing $g$.  Then we have the following inclusions of fixed-point sets:
\begin{equation} \label{e:fixedpointsets}
M^G \subset M^T \subset M^g \subset M^C.
\end{equation}

\begin{rem} \label{r:localtorsion}
If $T$ is a maximal torus in the compact connected Lie
group $G$ and $\dim T = \ell$, then
\[
H^*(BG) = H^*(BT)^{W_G} = \Q [u_1, \dots, u_{\ell}]^{W_G},
\]
where $W_G$ is the Weyl group of $T$ in $G$.
Thus, $H^*(BG)$ is an integral domain.  Let $Q$ be its field
of fractions.
For any $H^*(BG)$-module $V$, we define the localization of
$V$ with respect to the zero ideal in $H^*(BG)$ to be
\[
\hat{V} := V \otimes_{H^*(BG)} Q.
\]
It is easily verified that $V$ is $H^*(BG)$-torsion if and
only if $\hat{V} =0$.  For a $G$-manifold $M$, we call
$\hat{H}_G^*(M)$ the \emph{localized equivariant cohomology}
of $M$.
\end{rem}

\begin{lem}\label{l:gtor} Let $M$ be a $G$-manifold and $T$ a maximal torus of $G$.
If $H^*_T(M)$ is $H^*(BT)$-torsion, then
$H^*_G(M)$ is $H^*(BG)$-torsion.
\end{lem}
\begin{proof}
Recall that $H^*_G(M)$ is the subring of $H^*_T(M)$ consisting of the $W_G$-invariant elements.
Since $H^*_T(M)$ is $H^*(BT)$-torsion, there is $a\in H^*(BT)$ such
that $a\cdot 1_{H^*_T(M)}=0$. Consider the average of $a$ over the
Weyl group $W_G$ of $T$ in $G$,
$$
\tilde{a}=\frac{1}{|W_G|}(a+\omega_1 a+\cdots +\omega_r a) \in H^*(BG).
$$
Under  $\psi$, the element $\tilde{a}\cdot 1_{H^*_G(M)}$ goes to
$$
\frac{1}{|W_G|}(\omega_1 a+\cdots +\omega_r a) 1_{H^*_T(M)}.
$$
But $(\omega_j a)1_{H^*_T(M)}=\omega_j (a1_{H^*_T(M)})=0$ for any $j$.
Thus $\tilde{a}\cdot 1_{H^*_G(M)}=0$ in $H^*_G(M)$.
\end{proof}

\begin{prop}
\label{p:nfptorsion}
Let $G$ be a compact connected Lie group acting on a compact manifold $M$, and let $C$ be a conjugacy class in $G$.  If $U\subset M-M^C$ is an open G-subset, 
then the equivariant cohomology $H_G^*(U)$ is $H^*(BG)$-torsion.
\end{prop}

\begin{proof}
It follows from (\eqref{e:fixedpointsets}) that $U\subset M-M^C \subset M-M^T$.  
Since the inclusion map $U \to M-M^T$ is $T$-equivariant, and
$H_T^*(M-M^T)$ is $H^*(BT)$-torsion by
(\cite{guillemin-sternberg}, Th. 11.4.1),
$H_T^*(U)$ is also $H^*(BT)$-torsion.
By Lemma 1.2, $H_G^*(U)$ is $H^*(BG)$-torsion.
\end{proof}

In the rest of this section, ``torsion'' will mean
$H^*(BG)$-torsion.

\begin{thm}[Borel-type localization formula for a conjugacy
class] \label{t:borelconj}
Let $G$ be a compact connected Lie group acting on a compact
manifold $M$, and $C$ a conjugacy class in $G$.  Then the
inclusion $i: M^C \to M$ %of the fixed point set of the conjugacy class $C$
induces an isomorphism in localized
equivariant cohomology
\[
i^*: \hat{H}_G^*(M) \to \hat{H}_G^*(M^C).
\]
\end{thm}

\begin{proof}
Let $U$ be a $G$-invariant tubular neighborhood of $M^C$.
Then $\{U, M-M^C\}$ is a $G$-invariant open cover of $M$.
Moreover, ${H}_G^*(U) \simeq {H}_G^*(M^C)$ because
$U$ has the $G$-homotopy type of $M^C$.

By Prop.~\ref{p:nfptorsion}, $H_G^*(M-M^C)$
and $H^*_G(U\cap (M-M^C))$ are torsion.
Then in the localized equivariant Mayer-Vietoris sequence
\begin{align*}
\dots &\to \hat{H}_G^{*-1}(U\cap (M-M^C))\\
\to \hat{H}_G^*(M)
\to \hat{H}_G^*(M-M^C) \oplus \hat{H}_G^*(U) &\to
\hat{H}_G^{*}(U\cap (M-M^C)) \to \dots ,
\end{align*}
all the terms except $\hat{H}_G^*(M)$ and $\hat{H}_G^*(U)$
are zero.
It follows that
\[
\hat{H}_G^*(M) \to \hat{H}_G^*(U)\simeq \hat{H}^*_G(M^C)
\]
is an isomorphism of $H^*(BG)$-modules.
\end{proof}

When the group is a torus $T$, a conjugacy  class $C$ consist of a single element $t\in T$. If
$t$ is generator, then the fixed point set of $t$ is the same as the fixed
point set of the whole group $T$: $M^C=M^t=M^T$. 
In this case $M^C$ is smooth.
Thus Borel's localization theorem
follows from Theorem \ref{t:borelconj} by taking the conjugacy class $C=\{t\}$ in $T$.

\section{The equivariant Euler class}
\setcounter{thm}{0}

Suppose a compact connected Lie group $G$ acts on a smooth compact manifold
$M$. Let $C$ be a conjugacy class in $G$, and $M^C$ as before.
From now on we assume that $M^C$ is smooth with oriented normal bundle.
Denote by $i: M^C \to M$ the inclusion map and by $e_M\in H^*_G(M^C)$ the
equivariant Euler class of the normal bundle of $M^C$ in
$M$.

\begin{prop}\label{p:inverse}
Let $M$ be a compact connected oriented $G$-manifold.
Then the equivariant Euler class  $e_M$ of the normal bundle
of $M^C$ in $M$ is invertible in $\hat{H}^*_G(M^C)$.
\end{prop}

\begin{proof}
Fix a $G$-invariant Riemannian metric on $M$. Then the normal bundle
$\nu\to M^C$ is a $G$-equivariant vector bundle.
Let $\nu_0$ be the normal bundle minus the zero section.
Since $\nu_0$ is equivariantly diffeomorphic to an open
set in $M-M^C$, $\hat{H}^*_G(\nu_0)$ vanishes by Prop. \ref{p:nfptorsion}.
From the
Gysin long exact sequence in localized equivariant cohomology
$$
\cdots\to \hat{H}^*_G(\nu_0)\to \hat{H}^*_G(M^C)
\stackrel{\times e_M}\longrightarrow \hat{H}_G^*(M^C)
\to \hat{H}^*_G(\nu_0)\to\cdots
$$
it follows that multiplication
by the equivariant Euler class gives an automorphism of $\hat{H}^*_G(M^C)$.
Thus $e_M$ has an inverse in the ring $\hat{H}_G^*(M^C)$.
\end{proof}

Recall that the inclusion map $i:M^C\to M$ satisfies the identity
$$
i^*i_*(x)=x e_M, \qquad x \in H_G^*(M).
$$
in equivariant cohomology. 
In the localized equivariant cohomology $\hat{H}_G^*(M^C)$,
\[
i^*i_* \frac{i^*x}{e_M} = \frac{i^*x}{e_M} e_M = i^*x.\]
By Theorem \ref{t:borelconj}, $i^*$ is an isomorphism.  Hence,
\begin{eqnarray}\label{eqn:inverse}
i_*\left( \frac{i^*a}{e_M}  \right) = a.
\end{eqnarray}
for $a\in \hat{H}_G^*(M)$.

\section{Relative localization formula}
\setcounter{thm}{0}

Let $N$ be a $G$-manifold, $e_N$ the
equivariant Euler class of the normal bundle of $N^C$, and
$f: M \to N$  a $G$-equivariant map.
There is a
commutative diagram of maps
\begin{equation}\label{d:fixedpoints}
\square[M^C`M`N^C`N,; i_M`f^C`f`i_N]
\end{equation}
where $i_M$ and $i_N$ are inclusion maps and $f^C$ is the restriction of $f$ to $M^C$.
Let
\[
(f_G)_*: \hat{H}_G^*(M) \to \hat{H}_G^*(N), \qquad f_*^C:
\hat{H}_G^*(M^C) \to \hat{H}_G^*(N^C)
\]
be the push-forward maps in localized equivariant
cohomology.

\begin{thm}[Relative localization formula] \label{t:relativeconj}
Let $M$ and $N$ be compact oriented manifolds on which a compact connected Lie group $G$ acts, and
$f:M \to N$ a $G$-equivariant map.
For $a \in {H}_G^*(M)$,
\[
(f_G)_*a =(i_N^*)\inv f_*^C  \left( \frac{(f^C)^*e_N}{e_M} i_M^*
a\right)
\]
where the push-forward and restriction maps are in
localized equivariant cohomology.
\end{thm}

\begin{proof}
The commutative diagram (\ref{d:fixedpoints}),
induces a commutative diagram in localized equivariant
cohomology
\begin{equation} \label{d:pushforward}
\bfig
\square<1`1`1`1;800`400>[\hat{H}_G^*(M^C)`\hat{H}_G^*(M)`\hat{H}_G^*(N^C)`%
\hat{H}_G^*(N).;i_{M*}`f^C_*`(f_G)_*`i_{N*}]
\efig
\end{equation}

By eq. (\ref{eqn:inverse}) and the commutativity of the
diagram \eqref{d:pushforward},
\begin{align*}
(f_G)_*a &= (f_G)_* i_{M*} \left(\frac{1}{e_M} i_M^*   a\right) \\
&=i_{N*}f^C_* \left( \frac{1}{e_M} i_M^* a\right).
\end{align*}
Hence,
\begin{alignat*}{2}
i_N^* (f_G)_* a &= i_N^* i_{N*} f^C_* \left( \frac{1}{e_M} i_M^*
a\right) \\
&= e_N f^C_* \left( \frac{1}{e_M} i_M^* a\right) &&
  \\ %\qquad\text{(by \eqref{e:pushpull})} \\
&= (f^C)_* \left( \frac{(f^C)^*e_N}{e_M} i_M^* a\right) &&
\qquad\text{(projection formula)}.
\end{alignat*}
By Theorem \ref{t:borelconj}, $i_N^*$ is an isomorphism in localized equivariant
cohomology,
\[
(f_G)_*a = (i_N^*)\inv (f^C)_* \left( \frac{(f^C)^*e_N}{e_M} i_M^*
a\right).
\]
\end{proof}

If in Theorem \ref{t:relativeconj} we take the group $G$ to be a torus $T$ and the conjugacy class $C$ to be the conjugacy class of a generator $t$ for $T$, then $M^C = M^t = M^T$ and Theorem
\ref{t:relativeconj} specializes to the following formula of Lian, Liu, and Yau  \cite{lian-liu-yau}.

\begin{cor}[Relative localization formula for a torus
action] \label{c:relativetorus}
Let $M$ and $N$ be manifolds on which a torus $T$ acts, and
$f:M \to N$ a $T$-equivariant map with compact oriented
fibers. For $a \in \hat{H}_T^*(M)$,
\[
(f_T)_* a = (i_N^*)\inv (f^T)_* \left(
\frac{(f^T)^*e_N}{e_M} i_M^*a\right),
\]
where the push-forward and restriction maps are in
localized equivariant cohomology.
\end{cor}

When $N$ is a single point, Cor. \ref{c:relativetorus} reduces
to the Atiyah-Bott-Berline-Vergne localization formula.

\section{Applications to the Gysin homomorphism in ordinary cohomology}
\setcounter{thm}{0}

Let $G$ be a compact connected Lie group acting on a manifold $M$.
Denote by $M_G$ the homotopy quotient of $M$ by $G$,
and by $M^G$ the fixed point set of the action of $G$ on $M$.
Let $h_M: M \to M_G$ be the
inclusion of $M$ as a fiber of the bundle $M_G \to BG$ and
$i_M: M^G \to M$ the inclusion of the fixed point set $M^G$
in $M$.
The map $h_M$ induces a homomorphism in cohomology
\[
h_M^*: H_G^*(M) \to H^*(M).
\]
The inclusion $i_M$ induces a homomorphism in equivariant
cohomology
\[
i_M^*: H_G^*(M) \to H_G^*(M^G).
\]

A cohomology class $a \in H^*(M)$ is said to have an
\emph{equivariant extension} $\tilde{a} \in H_G^*(M)$ under
the $G$ action if under the restriction map
$h_M^*: H_G^*(M) \to H^*(M)$, the equivariant class
$\tilde{a}$ restricts to $a$.

Suppose $f: M\to N$ is a $G$-equivariant map of compact
oriented $G$-manifolds.  In this section we show that if a
class in $H^*(M)$ has an equivariant extension, then its
image under the Gysin map $f_*: H^*(M) \to H^*(N)$ in
ordinary cohomology can be computed from the relative
localization formulas (Cor.~\ref{c:relativetorus} or
Th.~\ref{t:relativeconj}).

We consider first the case of an action by a torus $T$.
Let $f_T: M_T \to N_T$ be the induced map of homotopy
quotients and $f^T: M^T\to N^T$ the induced map of fixed
point sets.  As before, $e_M$ denotes the equivariant Euler
class of the normal bundle of the fixed point set $M^T$ in $M$.

\begin{prop} \label{p:gysin}
Let $f: M\to N$ be a $T$-equivariant map of compact
oriented $T$-manifolds.  If a cohomology class $a \in
H^*(M)$ has an equivariant extension $\tilde{a} \in
H_T^*(M)$, then its image under the Gysin map $f_*: H^*(M)
\to H^*(N)$ is,
\begin{enumerate}
\item[1)] in terms of equivariant integration over $M$:
\[
f_*a = h_N^*f_{T*} \tilde{a},
\]
\item[2)] in terms of equivariant integration over the fixed
point set $M^T$:
\[
f_* a = h_N^* (i_N^*)\inv (f^T)_*
\left(
\frac{(f^T)^*e_N}{e_M} i_M^*\tilde{a}\right).
\]
\end{enumerate}
\end{prop}

\begin{proof}
The inclusions $h_M: M \to M_T$ and $h_N: N\to N_T$ fit
into a commutative diagram
\[
\square[M`M_T`N`N_T.;h_M`f`f_T`h_N]
\]
This diagram is Cartesian in the sense that $M$ is the
inverse image of $N$ under $f_T$.  Hence, the push-pull
formula $f_*h_M^*= h_N^*f_{T*}$ holds.  Then
\[
f_* a = f_*h_M^*\tilde{a} = h_N^* f_{T*} \tilde{a}.
\]

2) follows from 1) and the relative localization formula for
a torus action (Cor.~\ref{c:relativetorus}).
\end{proof}

Using the relative localization formula for a conjugacy
class, one obtains analogously a push-forward formula in
terms of the fixed point sets of a conjugacy class.
Now $h_M$ and $i_M$ are the inclusion maps
\[
h_M: M \to M_G, \qquad i_M: M^C \to M,
\]
$e_M$ is the equivariant Euler class of the normal bundle
of $M^C$ in $M$, and $f^C : M^C \to N^C$ is the induced map
on the fixed point sets of the conjugacy class $C$.

\begin{prop}
Let $f: M \to N$ be a $G$-equivariant map of compact
oriented $G$-manifolds.  Assume that the fixed point sets
$M^C$ and $N^C$ are smooth with oriented normal bundle.
For a class $a \in H^*(M)$ that has
an equivariant extension $\tilde{a} \in H_G^*(M)$,
\[
f_* a = h_N^* (i_N^*)\inv(f^C)_* \left(
\frac{(f^C)^*e_N}{e_M} i_M^*\tilde{a}\right).
\]
\end{prop}

\section{Example: the Gysin homomorphism of flag manifolds}
\setcounter{thm}{0}

Let $G$ be a compact connected Lie group with maximal torus
$T$, and $H$ a closed subgroup of $G$ containing $T$.
In \cite{akyildiz-carrell} Akyildiz and Carrell compute
the Gysin homomorphism for the canonical projection $f:
G/T \to G/H$.  In this section we deduce the formula of
Akyildiz and Carrell from the relative localization formula
in equivariant cohomology.

Let $N_G(T)$ be the normalizer of the torus $T$ in the
group $G$.  The Weyl group $W_G$ of $T$ in $G$ is
$W_G=N_G(T)/T$.
We use the same letter $w$ to denote an element of the Weyl
group $W_G$ and a lift of the element to the normalizer
$N_G(T)$.  The Weyl group $W_G$ acts on $G/T$ by
\[
(gT)w= gwT \qquad\text{for $gT \in G/T$ and\ } w\in W_G.
\]
This induces an action of $W_G$ on the cohomology ring
$H^*(G/T)$.

We may also consider the Weyl group $W_H$ of $T$ in $H$.
By restriction the Weyl group $W_H$ acts on $G/T$ and on
$H^*(G/T)$.

To each character $\gc$ of $T$ with representation space
$\C_{\gc}$, one associates a complex line bundle
\[
L_{\gc} := G \times_T \C_{\gc}
\]
over $G/T$.  Fix a set $\tri^+(H)$ of positive roots for
$T$ in $H$, and extend $\tri^+(H)$ to a set $\tri^+$ of
positive roots for $T$ in $G$.

\begin{thm}[\cite{akyildiz-carrell}]
\label{t:akyildiz-carrell}
The Gysin homomorphism $f_*: H^*(G/T) \to H^*(G/H)$ is
given by, for $a \in H^*(G/T)$,
\[
f_* a =
\dfrac{\sum_{w\in W_H} (-1)^w w\cdot a}
{\prod_{\ga \in \tri^+(H)} c_1(L_{\ga})}.
\]
\end{thm}

\begin{rem}
There are two other ways to obtain this formula.  First,
using representation theory, Brion \cite{brion96} proves a
push-forward formula for flag bundles that includes
Th.~\ref{t:akyildiz-carrell} as a special case.  Secondly,
since $G/T \to G/H$ is a fiber bundle with equivariantly
formal fibers, the method of \cite{tu-gysin} using the ABBV
localization theorem also applies.
\end{rem}

To deduce Th.~\ref{t:akyildiz-carrell} from
Prop.~\ref{p:gysin} we need to recall a few facts about the
cohomology and equivariant cohomology of $G/T$ and $G/H$
(see \cite{tu-gysin}).

\subsection*{Cohomology ring of $BT$}

Let $ET \to BT$ be the universal principal $T$-bundle.  To
each character $\gc$ of $T$, one associates a complex line
bundle $S_{\gc}$ over $BT$ and a complex line bundle $L_{\gc}$ over $G/T$:
\[
S_{\gc} := ET \times_T \C_{\gc}, \qquad
L_{\gc} := G \times_T  \C_{\gc} .
\]
For definiteness, fix a basis $\chi_1, \dots, \chi_{\ell}$
for the character group $\hat{T}$, where we write the
characters additively, and set
\[
u_i = c_1(S_{\chi_i}) \in H^2(BT), \qquad z_i =
c_1(L_{\chi_i}) \in H^2 (G/T).
\]
Let $R= \Sym(\hat{T})$ be the symmetric algebra over $\Q$
generated by $\hat{T}$.  The map $\gc \mapsto c_1(S_{\gc})$
induces an isomorphism
\[
R= \Sym(\hat{T}) \overset{~}{\to} H^*(BT) = \Q [ u_1, \dots,
u_{\ell}].\]
The map $\gc \mapsto c_1(L_{\gc})$ induces an isomorphism
\[
R = \Sym(\hat{T}) \overset{~}{\to} \Q [ z_1, \dots, z_{\ell}].
\]
The Weyl groups $W_G$ and $W_H$ act on the characters of $T$
and hence on $R$: for $w\in W_G$ and $\gc \in \hat{T}$,
\[
(w\cdot \gc)(t) = \gc(w\inv t w).
\]

\subsection*{Cohomology rings of flag manifolds}

The cohomology rings of $G/T$ and $G/H$ are described in
\cite{borel53}:
\begin{align*}
H^*(G/T) &\simeq \frac{R}{(R_+^{W_G})} \simeq
\frac{\Q[z_1, \dots, z_{\ell}]}{(R_+^{W_G})},\\
H^*(G/H) & \simeq \frac{R^{W_H}}{(R_+^{W_G})} \simeq
\frac{\Q[z_1, \dots, z_{\ell}]^{W_H}}{(R_+^{W_G})},
\end{align*}
where $(R_+^{W_G})$ denotes the ideal generated by the
$W_G$-invariant homogeneous polynomials of positive degree.

The torus $T$ acts on $G/T$ and $G/H$ by left
multiplication.  
For each character $\chi$ of $T$, let $K_{\chi}:=(L_{\chi})_T$
be the homotopy quotient of the bundle $L_{\chi}$ by the torus $T$.
Then $K_{\chi}$ is a complex line bundle over $(G/T)_T$.
Their equivariant cohomology rings are
(see %\cite{brion}, 
\cite{tu-gysin})
\begin{align*}
H_T^*(G/T) &= \frac{ \Q [u_1, \dots, u_{\ell}, y_1, \dots,
y_{\ell}]}{J},\\
H_T^*(G/H) &= \frac{ \Q [u_1, \dots, u_{\ell}]\otimes(\Q[y_1, \dots,
y_{\ell}]^{W_H})}{J},\\
\end{align*}
where $y_i = c_1(K_{\chi_i}) \in H_T^*(G/T)$ 
and $J$ denotes the ideal generated by $q(y) - q(u)$ for
$q \in R_+^{W_G}$.

\subsection*{Fixed point sets}

The fixed point sets of the $T$-action on $G/T$ and on
$G/H$ are the Weyl group $W_G$ and the coset space
$W_G/W_H$ respectively.
Since these are finite sets of points,
\begin{gather*}
H_T^*(W_G) = \oplus_{w\in W_G} H_T^*(\{ w\}) \simeq
\oplus_{w\in W_G} R, \\
H_T^*(W_G/W_H) = \oplus_{w W_H \in W_G/W_H} R.
\end{gather*}
Thus, we may view an element of $H_T^*(W_G)$ as a function
from $W_G$ to $R$, and an element of $H_T^*(W_G/W_H)$ as a
function from $W_G/W_H$ to $R$.

Let $h_M: M \to M_T$ be the inclusion of $M$ as a fiber in
the fiber bundle $M_T \to BT$ and $i_M: M^T \to M$ the
inclusion of the fixed point set $M^T$ in $M$.
Note that $i_M$ is $T$-equivariant and induces a
homomorphism in $T$-equivariant cohomology,
$i_M^*: H_T^*(M) \to H_T^*(M^T)$.
In order to apply Prop. \ref{p:gysin}, we need to know how to
calculate the restriction maps
\[
h_M^*: H_T^*(M) \to H^*(M) \quad\text{and}\quad i_M^*:
H_T^*(M) \to H_T^*(M^T)
\]
as well as the equivariant Euler class $e_M$ of the normal
bundle to the fixed point set $M^T$, for $M= G/T$ and
$G/H$.  This is done in \cite{tu-gysin}.

\subsection*{Restriction and equivariant Euler class
formulas for $G/T$}

Since $h_M^*: H_T^*(M) \to H^*(M)$ is the restriction to a
fiber of the bundle $M_T \to BT$, and the bundle
$K_{\chi_i}= (L_{\chi_i})_T$ on $M_T$ pulls back to
$L_{\chi_i}$ on $M$,
\begin{equation} \label{e:restricth}
h_M^*(u_i) = 0, \qquad h_M^*(y_i) =
h_M^*(c_1(K_{\chi_i})) = c_1(L_{\chi_i}) = z_i.
\end{equation}

Let $i_w: \{ w\} \to G/T$ be the inclusion of the fixed
point $w\in W_G$ and
\[
i_w^*: H_T^*(G/T) \to H_T^*(\{w\}) = R
\]
the induced map in equivariant cohomology.
By (\cite{tu-gysin}), for $p(y) \in H_T^*(G/T)$,
\begin{equation} \label{e:restricti}
i_w^* u_i = u_i, \quad i_w^* p(y) = w\cdot p(u), \quad i_w^*
c_1(K_{\gc}) = w\cdot c_1(S_{\gc}).
\end{equation}
Thus, the restriction of $p(y)$ to the fixed point set
$W_G$ is the function $i_M^*p(y): W_G \to R$ whose value
at $w\in W_G$ is
\begin{equation} \label{e:restrictim}
(i_M^*p(y))(w)= w \cdot p(u).
\end{equation}

The equivariant Euler class of the normal bundle to the
fixed point set $W_G$ assigns to each $w\in W_G$ the
equivariant Euler class of the normal bundle $\nu_w$ at
$w$; thus, it is also a function $e_M: W_G \to R$.
By (\cite{tu-gysin}),
\begin{equation} \label{e:eulerm}
e_M(w) = e^T(\nu_w) = w\left( \prod_{\ga \in \tri^+}
c_1(S_{\ga})\right)
= (-1)^w \prod_{\ga \in \tri^+} c_1 (S_{\ga}).
\end{equation}

\subsection*{Restriction and equivariant Euler class
formulas for $G/H$}

For the manifold $M= G/H$, the formulas for the restriction
maps $h_N^*$ and $i_N^*$ are the same as in
\eqref{e:restricth} and \eqref{e:restricti}, except that now
the polynomial $p(y)$ must be $W_H$-invariant.  In
particular,
\begin{equation} \label{e:restricthn}
h_N^*(u_i) = 0, \quad h_N^* p(y) = p(z), \quad
h_N^*(c_1(K_{\gc})) = c_1(L_{\gc}),
\end{equation}
and
\begin{equation} \label{e:restrictin}
(i_N^* p(y))(w W_H) = w \cdot p(u).
\end{equation}

If $\gc_1, \dots, \gc_m$ are characters of $T$ such that
$p(c_1(K_{\gc_1}), \dots, c_1(K_{\gc_m}))$ is invariant
under the Weyl group $W_H$, then
\begin{equation} \label{e:restrictin2}
(i_N^*p (c_1(K_{\gc_1}), \dots, c_1(K_{\gc_m})))(w W_H) =
w\cdot p(c_1(S_{\gc_1}), \dots, c_1(S_{\gc_m})).
\end{equation}

The equivariant Euler class of the normal bundle of the
fixed point set $W_G/W_H$ is the function
$e_N: W_G/W_H \to R$ given by
\begin{equation} \label{e:eulern}
e_N(w W_H) = w \cdot \left(
\prod_{\ga \in \tri^+ - \tri^+(H)} c_1(S_{\ga})\right).
\end{equation}

\begin{proof}[Proof of Th.~\ref{t:akyildiz-carrell}]
With $M=G/T$ and $N= G/H$ in Prop.~\ref{p:gysin}, let
\[
p(z) \in H^*(G/T) = \Q[z_1, \dots,
z_{\ell}]/(R_+^{W_G}).
\]
It is the image of
$p(y) \in H_T^*(G/T)$ under the restriction map
$h_M^*: H_T^*(G/T) \to H^*(G/T)$.
By Prop. \ref{p:gysin},
\begin{equation} \label{e:pushforward}
f_* p(z) = f_* h_M^* p(y) = h_N^*f_{T*} p(y)
\end{equation}
and
\[
f_{T*} p(y) = (i_N^*)\inv (f^T)_*
\left( \frac{(f^T)^*e_N}{e_M} i_M^*
p(y) \right).
\]

By Eq.~\eqref{e:restrictim}, \eqref{e:eulerm},
and \eqref{e:eulern},
for $w \in W_G$,
\[
(i_M^*p(y))(w) = i_w^* p(y) = w\cdot p(u),
\]
and
\begin{align*}
\left( \frac{(f^T)^*e_N}{e_M} \right)(w) &=
\frac{e_N(w W_H)}{e_M(w)} =
w\cdot
\left(
\frac{\prod_{\ga \in \tri^+ - \tri^+(H)} c_1 (S_{\ga})}
{\prod_{\ga \in \tri^+} c_1( S_{\ga})} \right) \\
&= \frac{1}{w\cdot \left(
\prod_{\ga \in \tri^+(H)} c_1(S_{\ga}) \right) }.
\end{align*}

To simplify the notation, define temporarily the function
$k: W_G \to R$ by
\[
k(w) = w\cdot \left(
\frac{p(u)}{\prod_{\ga \in \tri^+(H)} c_1(S_{\ga})} \right).
\]
Then
\begin{equation} \label{e:k1}
f_{T*} p(y) = (i_N^*)\inv (f^T)_*(k).
\end{equation}

Now $(f^T)_*(k) \in H_T^*(W_G/W_H)$ is the function$:
W_G/W_H \to R$ whose value at the point $w W_H$ is obtained
by summing $k$ over the fiber of $f^T: W_G \to W_G/W_H$ above
$w W_H$.  Hence,
\begin{align*}
((f^T)_* k) (w W_H) &= \sum_{w v \in w W_H} wv \cdot \left(
\frac{p(u)}{\prod_{\ga \in \tri^+(H)} c_1(S_{\ga})}
\right) \\
&= w\cdot \sum_{v\in W_H} v\cdot \left(
\frac{p(u)}{\prod_{\ga \in \tri^+(H)} c_1(S_{\ga})}
\right) .
\end{align*}
By \eqref{e:restrictin2}, the inverse image of this
expression under $i_N^*$ is
\begin{equation} \label{e:k2}
(i_N^*)\inv (f^T)_* k = \sum_{v\in W_H} v\cdot \left(
\frac{p(y)}{\prod_{\ga \in \tri^+(H)} c_1(K_{\ga})}
\right) .
\end{equation}
Finally, combining \eqref{e:pushforward}, \eqref{e:k1},
\eqref{e:k2} and \eqref{e:restricthn},
\[
f_*p(z) = h_N^*(f_T)_* p(y) = \sum_{v\in W_H} v\cdot
\left(
\frac{p(z)}{\prod_{\ga \in \tri^+(H)} c_1(L_{\ga})}
\right) .
\]
\end{proof}

\end{document}